\magnification1200
\parindent0pt
\parskip3pt
\baselineskip=1.5\baselineskip
\input pictex
\input amssym.tex
\def\dzien{\number\day\space%
\ifcase\month \or
 stycznia \or lutego \or marca \or kwietnia \or maja \or czerwca \or
 lipca \or sierpnia \or wrzeœnia \or paŸdziernika \or listopada \or
 grudnia \fi
 \number\year}
\newcount\t
\newcount\h
\newcount\m
\def\godzina{\t=\number\time \h=\t
\divide\h by60 \the\h:%
\multiply\h by60
\advance\t by-\h \the\t}


\def\title#1{\centerline{\bf#1}}
\def\section#1{\bigskip{\bf#1.}\ }

\def\frac#1#2{{#1\over#2}}
\def\binom#1#2{{#1\choose#2}}
\def\tw{\hat{\theta}_{w}} 
\def\tc{\hat{\theta}_{c}} 
\def\cU{{\cal U}}

\title{Cost Issue in Estimation of Proportion}
\title{in a Finite Population Divided Among Two Strata}
\bigskip
\centerline{Dominik Sieradzki, Wojciech Zieli\'nski}
\centerline{Department of Econometrics and Statistics}
\centerline{Warsaw University of Life Sciences}
\centerline{Nowoursynowska 159, 02-776 Warszawa, Poland}
\centerline{e-mail: dominik$\_$sieradzki@sggw.pl}
\centerline{e-mail: wojciech$\_$zielinski@sggw.pl}
\bigskip

\noindent{\bf Abstract. } The problem of estimation of the proportion of units with a given attribute in a~finite population is considered. From the population a sample is drawn due to the simple random sampling without replacement. There are limited funds for conducting survey sample. Suppose that the population is divided into two strata. The question now arises: how should sample sizes be chosen from each strata to obtain the best estimation of proportion without exceeding the budget planned. In the paper it is shown, that with the appropriate sample allocation the variance of the stratified estimator may be reduced up to $30\%$ off of the standard, unstratified estimator.
\bigskip


\bigskip

\section{\bf 1. Introduction}
Consider a population $\cU=\{u_1,\ldots,u_N\}$ which contains a finite number $N$ of units. Let $M$ denote an unknown number of objects in the population which have an interesting property. The aim of the study is to estimate $M$, or equivalently, the fraction $\theta={M\over N}$. The sample of size $n$ is drawn due to the simple random sampling without replacement. Let $\xi$ be a random variable describing a number of objects with the property in the sample. The random variable $\xi$ has a hypergeometric distribution (Johnson, Kotz 1969):
$$P_{\theta,N,n} \left\{ \xi = x \right\} = \frac{\binom{\theta N}{x} \binom{\left(1-\theta\right) N}{n-x}}{\binom{N}{n}},\eqno(1.1)$$
for integer $x$ from the set $\left\{\max\{0, n - \left(1-\theta\right) N\},\ldots, \min\{n, \theta N\}\right\} $. It is known that the unbiased estimator with minimal variance of the parameter $\theta$ is $\hat{\theta}_c = \frac{\xi}{n}$ and its variance equals
$D^2_{\theta}\tc = \frac{ \theta (1- \theta)}{ n} \frac{N-n}{N-1}.$

Suppose that the population $\cU$ is divided into two strata $\cU_1$ and $\cU_2$, $\cU= \cU_1 \cup \cU_2$, $\cU_1\cap\cU_2=\emptyset$, of $N_1$ and $N_2$ units, respectively. The problem lays in finding sample sizes $n_1$ and $n_2$ from the first and the second strata in order to obtain best estimation of $\theta$. Assume that the overall cost of sampling is limited to a given number $C$. In the second section of the paper a stratified estimator $\tw$ of $\theta$ is proposed and its variance is analyzed. This variance is compared with the variance of $\tc$. It is shown the existence of $n_1$ and $n_2$ such that the variance of $\tw$ is smaller than the variance of $\tc$ for all $\theta$. In the third section some numerical results are presented.

Results of the paper may be easily generalized to an arbitrary number of strata.


\section{\bf 2. Stratified estimator}
Let the population $\cU$ be divided into two strata $\cU_1$ and $\cU_2$. In each strata proportions $\theta_1$ and $\theta_2$ of objects with a given attribute are unknown. Let the contribution of the first strata be $w_1$, i.e. $w_1=N_1/N$. The overall proportion $\theta$ equals
$$
\theta = w_1 \theta_1 + w_2 \theta_2,\eqno(2.1)
$$
where $w_2 = 1 - w_1$. The estimate of $\theta$ is taken as (Horgan 2006):
$$
\tw=w_1\frac{\xi_1}{n_1}+w_2\frac{\xi_2}{n_2},\eqno(2.2)
$$
where $n_1$ and $n_2$ denote sample sizes from the first and second strata, respectively. Now, there are two random variables describing a number of units with a~particular attribute in samples drawn from each strata:
$$
\xi_1 \sim H\left(N_1,\theta_1 N_1, n_1 \right), \quad \xi_2 \sim H\left(N_2,\theta_2 N_2, n_2 \right).\eqno(2.3)
$$
The variance of the estimator $\tw$ equals:
$$
D^2_{\theta_1,\theta_2}\tw =\frac{w_1^2}{ n_1^2}  D^2_{\theta_1} \xi_1 +\frac{w_2^2}{ n_2^2}  D^2_{\theta_2} \xi_2 = w_1^2\frac{ \theta_1 (1- \theta_1)}{ n_1} \frac{N_1-n_1}{N_1-1} +w_2^2\frac{ \theta_2 (1- \theta_2)}{ n_2} \frac{N_2-n_2}{N_2-1}.\eqno(2.4) 
$$
Consider a cost of sampling. Suppose that the individual cost of sampling from the first strata is $c_1$ and from the second one is $c_2$. The cost of sampling equals $c_1 n_1 + c_2 n_2$ and is limited by a given value $C$, i.e.
$$
c_1 n_1 + c_2 n_2\leq C.\eqno(2.5)
$$
It is assumed that $c_1N_1+c_2N_2>C$. The whole sample size equals $n=n_1+n_2$. The question is: how should $n_1$ and $n_2$ be chosen in order to obtain the best estimate of $\theta$ with the restriction (2.5)? This problem is known as a sample allocation between strata. One of known solutions of that problem is optimal allocation (Cochran 1977, Neyman 1934). To minimize the variance of $\tw$ for given costs $c_1$ and $c_2$ and for given $\theta_1$ and $\theta_2$ sample sizes $n_1$ and $n_2$ are:
$$
n_k =n {N_k\sqrt{\theta_k (1-\theta_k)/c_k} \over \sum_{i=1}^2 N_i \sqrt{\theta_i (1-\theta_i)/c_i} }, \quad (k = 1,2), \hbox{ where }
n=C {\sum_{i=1}^2 N_i\sqrt{\theta_i (1-\theta_i)/c_i} \over \sum_{i=1}^2 N_i \sqrt{\theta_i (1-\theta_i)c_i} }.\eqno(2.6)
$$
Optimal allocation requires knowledge of the parameters $\theta_1$ and $\theta_2$. They are known accurately if the population was subjected to exhaustive sampling. Usually values $\theta_1$ and $\theta_2$ are estimated from a preliminary sample. In some cases fairly good estimates of $\theta_1$ and $\theta_2$ are available from past experience (Armitage 1947). Because of these difficulties values $n_1$ and $n_2$ could be far from optimal, hence stratified random sampling may have higher variance than $D^2_{\theta}\tc$ (Cochran 1977; Hansen, Hurwitz 1946).

Since the aim of the study is to estimate the overall fraction $\theta$, hence the parameter $\theta_1$ will be considered as a~nuisance one. This parameter will be eliminated by appropriate averaging. Note that for a given  $\theta \in [0,1]$, parameter $\theta_1$ is a fraction $M_1/N_1$ (it is treated as a number, not as a random variable) from the set (we assume that $w_1 \le w_2$)
$$
{\cal A}=\left\{a_\theta,a_\theta+ \frac{1}{N_1},a_\theta+ \frac{2}{N_1}, \ldots , b_\theta  \right\},\eqno(2.7)
$$
where
$$a_\theta = \max\left\{0,\frac{\theta-w_2}{w_1}\right\} \quad \hbox{and} \quad b_\theta= \min\left\{1,\frac{\theta}{w_1}\right\}\eqno(2.8)
$$
and let $L_\theta$ be the cardinality of ${\cal A}$ (Zieli{\'n}ski 2016).

\bigskip \noindent
{\bf Theorem 1.} The estimator $\tw$ is an unbiased estimator of $\theta$ (Sieradzki, Zieli{\'n}ski 2017).

\bigskip \noindent
Averaged variance of the estimator $\tw$ equals:
$$\eqalign{
D_\theta^2 \tw& =D_\theta^2 \left( w_1 {\xi_1\over n_1} + w_2 {\xi_2\over n_2} \right)=  \cr
&= {1\over L_\theta} \sum_{\theta_1\in{\cal A}} \left(  \left({w_1\over n_1}\right)^2 D_{\theta_1}^2\xi_1  +   \left({w_2\over n_2}\right)^2 D_{{\theta-w_1 \theta_1\over w_2}}^2\xi_2  \right)     = \cr
&= {1\over L_\theta} \sum_{\theta_1\in{\cal A}} \left[  {w_1^2\over n_1} \theta_1(1-\theta_1) {N_1-n_1\over N_1-1}  +   {w_2^2\over n_2}   {\theta-w_1 \theta_1\over w_2} \left(1-{\theta-w_1 \theta_1\over w_2}\right) {N_2-n_2\over N_2-1}  \right] . \cr
}\eqno(2.9)
$$
Detailed analysis of the variance $D^2_{\theta} \tw$ is given in Sieradzki, Zieli{\'n}ski (2017). Afterwards, derivation of sample allocation between strata is considered. There is a need to find such values of $\left(n_1^{opt}, n_2^{opt}\right)$ which minimize $\max_\theta D^2_{\theta} \tw$ and $c_1 n_1^{opt} + c_2 n_2^{opt} \le C$.
The variance of $\tw$, having regard to the cost is as follows:
$$\eqalign{
D_\theta^2 \tw &=
\frac{1}{L_\theta} \sum_{\theta_1\in{\cal A}}\left[\frac{w_1^2}{n_1} \theta_1(1-\theta_1) \frac{N_1-n_1}{N_1-1} \right.\cr
&\left.+\frac{w_2^2}{(C-c_1n_1)/c_2}\frac{\theta-w_1 \theta_1}{w_2}\left(1-\frac{\theta-w_1 \theta_1}{w_2}\right)\frac{N_2-(C-c_1n_1)/c_2}{N_2-1}  \right] . \cr
}\eqno(2.10)
$$
\noindent
For $0<\theta<w_1$ we have (here $n_2=(C-c_1n_1)/c_2$)
$$\eqalign{
D_\theta^2 \tw =& \frac{\theta}{6N} \left[\frac{(3N_1-1)(N_1-n_1)}{n_1(N_1-1)}+\frac{(3N_2-1)(N_2-n_2)}{n_2(N_2-1)}- \right.\cr
& \left. 2N\left(\frac{N_1-n_1}{n_1(N_1-1)}+\frac{N_2-n_2}{n_2(N_2-1)}\right)\theta\right] . \cr}\eqno(2.11)$$
For $w_1<\theta<w_2$ we have
$$\eqalign{
D_\theta^2 \tw=& \frac{2 n_1 n_2 w_1 (N+1)+N_1 N_2 ((n_1+n_2) w_1-3n_1)-N_1 (n_2 w_1+n_1 (1-w_1))}{6 N n_1 n_2 (N_2-1)}
+\cr &\frac{N_2-n_2}{n_2}\frac{\theta(1-\theta)}{N_2-1}.\cr} \eqno(2.12)$$
To obtain explicit formula for the variance of $\tw$ for $1-w_1<\theta<1$ it is enough to replace $\theta$ by $1-\theta$ in $(2.11)$.

Depending on $n_1$ $\max_\theta D^2_{\theta} \tw$ is achieved at $1/2$ or $\theta^* \in (0,w_1)$. To prove that, it is enough to find $\theta^*$. Let $\widetilde{\theta}$ denote $
 \theta$ which maximizes $D^2_{\theta} \tw$, i.e. $\widetilde{\theta} = 1/2$ or $\widetilde{\theta} = \theta^*
 $   Since $D_\theta^2\tw$ is a quadratic function of $\theta$ hence after some elementary calculations we obtain the formula for $\theta^*$
$$\theta^*=\frac{1}{4N}\frac{\frac{(3N_1-1)(N_1-n_1)}{n_1(N_1-1)}+\frac{(3N_2-1)(N_2-n_2)}{n_2(N_2-1)}}{\frac{N_1-n_1}{n_1(N_1-1)}+\frac{N_2-n_2}{n_2(N_2-1)}}. \eqno(2.13)$$
The maximal value of the variance of $\tw$ equals
$$D_{\theta^*}^2\tw=\frac{1}{48N^2}\frac{\left(\frac{(3N_1-1)(N_1-n_1)}{n_1(N_1-1)}+\frac{(3N_2-1)(N_2-n_2)}{n_2(N_2-1)}\right)^2}{\frac{N_1-n_1}{n_1(N_1-1)}+\frac{N_2-n_2}{n_2(N_2-1)}} \eqno(2.14)$$
or
$$D_{0.5}^2\tw=\frac{1}{6N^2}\left(N_1\frac{N_1-n_1}{n_1}+\frac{3N_2^2-(N_1+1)^2+1}{2N_1}\frac{N_2-n_2}{n_2}\right) \eqno(2.15)$$
depending on $n_1$.

\bigskip
{\bf Theorem 2.} Maximal variance $D_\theta^2\tw$ equals
$$\eqalign{
&D_{\theta^*}^2\tw \hbox{ for } w_1>w_1^*\cr
\noalign{\hbox{and}}
&\cases{
D_{\theta^*}^2\tw&for $0\leq n_1\leq n_1^*$\cr
D_{0.5}^2\tw&for $n_1^*\leq n_1\leq\min\{C/c_1,N_1\}$\cr
}\hbox{ for } w_1\leq w_1^*, \cr} \eqno(2.16)
$$
where $w_1^*\in[0,0.5]$ is the solution of $\lim_{n_1\to\frac{C}{c_1}}\frac{D_{0.5}^2\tw}{D_{\theta^*}^2\tw}=1$ and $n_1^*\in[0,\min\{C/c_1,N_1\}]$ is the solution of $D_{\theta^*}^2\tw=D_{0.5}^2\tw$.

\medskip
{\it Proof.} Assume that $n_1$ is a continuous variable. Only $n_1^*\in[0,\min\{C/c_1,N_1\}]$ is considered.

The derivative of $D_{0.5}^2\tw$ with respect to $n_1$ is proportional to $$-\frac{N_1^2}{n_1^2}+\frac{\Lambda}{(C-c_1n_1)^2},\eqno(2.17)$$ where $\Lambda$ is a positive constant. The derivative is negative for small $n_1$ and is positive for large $n_1$. Hence the variance $D_{0.5}^2\tw$ is a bathtub-shaped function of $n_1$.

In a similar way it may be shown that $D_{\theta^*}^2\tw$ is a bathtub-shaped.

Since
$$\lim_{n_1\to 0}\frac{D_{0.5}^2\tw}{D_{\theta^*}^2\tw}=\frac{8Nw_1(Nw_1-1)}{(3Nw_1-1)^2}<1, \hbox{ for all }w_1\in[0,0.5], \eqno(2.18)$$
and
$$\lim_{n_1\to\frac{C}{c_1}}\frac{D_{0.5}^2\tw}{D_{\theta^*}^2\tw}=\frac{4(N (1-w_1)-1) (N (3 - 6 w_1 + 2 w_1^2)-2 w_1)}{(3 N (1-w_1)-1)^2 w_1}=\cases{\geq1& for $w_1\leq w_1^*$,\cr <1& for $w_1>w_1^*$,\cr} \eqno(2.19)$$

hence we obtain the thesis. Exemplary values of $w_1^*$ are given in Table 1.
$$\vbox{\tabskip1em \offinterlineskip \halign{
\strut\hfil$#$\hfil&&#\vrule&\hfil$#$\hfil\cr
\multispan{5}{\bf Table 1.} Values of $w_1^{*}$\hfil\cr\noalign{\vskip5pt}
N&& w_1^{*}&&N_1  \cr \noalign{\hrule}
100\phantom{00000}&&0.463384&&46\phantom{00000}\cr
1000\phantom{0000}&&0.464030&&464\phantom{0000}\cr
10000\phantom{000}&&0.464094&&4640\phantom{000}\cr
100000\phantom{00}&&0.464101&&46410\phantom{00}\cr
1000000\phantom{0}&&0.464102&&464101\phantom{0}\cr
10000000&&0.464102&&4641016\cr
}}$$

\bigskip

To find optimal $n_1$ it is enough to minimize the maximal variance with respect to $n_1$.

\bigskip
{\bf Theorem 3.} For $w_1\leq w_1^*$ the optimal allocation of the sample is ($n_2^{opt}=\left(C-c_1n_1^{opt}\right)/c_2$):
$$
n_1^{opt}=\frac{C\sqrt{(N_2-1)} w_1}{c_1\sqrt{(N_2-1)} w_1+\sqrt{c_1 c_2w_2(N (w_1^2-3w_1+1.5)-w_1)}}, \eqno(2.20)$$

\medskip
{\it Proof.} Note that $D_{\theta^*}^2\tw$ is a decreasing function of $n_1$ for $0\leq n_1\leq n_1^*$. For $n_1=n_1^*$ we have
$$D_{\theta^*}^2\tw=D_{0.5}^2\tw . \eqno(2.21)$$

Since $\frac{N_1-n_1}{n_1}$ is decreasing in $n_1$ and $\frac{N_2-n_2}{n_2}$ is increasing in $n_1$, $D_{0.5}^2\tw$ has a minimum. To find an optimal $n_1$ it is enough to solve the equation (assuming that $n_1$ is a continuous variable)
$$\frac{\partial}{\partial n_1}D_{0.5}^2\tw=0. \eqno(2.22)$$

For $w_1 > w_1^{*}$ to find optimal $n_1$ it is necessary to solve the following equation
$$\frac{\partial}{\partial n_1}D_{\theta^*}^2\tw=0. \eqno(2.23)$$

The closed formula for optimal $n_1$ is available, but its form is very complicated and useless. A numerical solution is suggested.

\bigskip

To compare variances of the estimators $\tw$ and $\tc$, it is necessary to determine sample size for the estimator $\tc$. When simple random sampling from the whole population is applied, there is no information about which strata given object is drawn from. Hence the number of objects drawn the first strata is a random variable. Denote this random variable by $\eta_1$. Its distribution is a hypergeometric $H(N,w_1 N, n)$. The expected cost of the sample of size $n$ is
$$\sum_{k=0}^n\left(c_1k+c_2(n-k)\right)P\left\{\eta_1=k\right\}=(w_1c_1+(1-w_1)c_2)n \eqno(2.24
)$$
and the expected sample size $n_c$ for the estimator $\tc$ is
$$n_c = \frac{C}{w_1 c_1 + (1-w_1)c_2} . \eqno(2.25)$$

\bigskip
{\bf Theorem.} For $n_1^{opt}$ and $n_c= \frac{C}{w_1 c_1 + (1-w_1)c_2}$ $$D_{0.5}^2\tw\leq D_{0.5}^2\tc\hbox{ for }w_1\leq w_1^*\quad\hbox{and}\quad D_{\theta^*}^2\tw\leq D_{0.5}^2\tc\hbox{ for }w_1> w_1^* \eqno(2.26)$$

\medskip
{\it Proof.} The maximal variance of the estimator $\tc$ equals
$$D_{0.5}^2\tc=\frac{c_1N_1 + c_2N_2}{4C(N-1)}-\frac{1}{4(N-1)}. \eqno(2.27)$$
If $w_1\leq w_1^*$ then the maximal variance of the estimator $\tw$ equals
$$\eqalign{
D_{0.5}^2\tw=&
\frac{\left(w_1\sqrt{c_1 N (N_2-1)}+\sqrt{c_2 N_2 (N (w_1^2-3 w_1+1.5)-w_1)}\right)^2}{6 C N (N_2-1)}-\cr
&\hskip2em\frac{(1.5N-2 N_1-2 w_1)}{6N (N_2-1)}\cr
} \eqno(2.28)$$

We have:
$$\eqalign{
&\frac{\left(w_1\sqrt{c_1 N (N_2-1)}+\sqrt{c_2 N_2 (N (w_1^2-3 w_1+1.5)-w_1)}\right)^2}{6 C N (N_2-1)}-\frac{(1.5N-2 N_1-2 w_1)}{6N (N_2-1)}\leq\cr
&\frac{\left(w_1\sqrt{c_1 N N_2}+\sqrt{c_2 N_2 N (w_1^2-3 w_1+1.5)}\right)^2}{6 C N (N_2-1)}-\frac{(1.5N-2 N_1)}{6N (N_2-1)}\leq\cr
&\frac{N_2}{6 C  (N_2-1)}\left(w_1\sqrt{c_1}+\sqrt{c_2(w_1^2-3 w_1+1.5)}\right)^2+\frac{w_1}{3(N-1)}-\frac{1}{4(N-1)} . \cr
} \eqno(2.29)$$
Now it is enough to show that
$$\frac{N_2}{6 C  (N_2-1)}\left(w_1\sqrt{c_1}+\sqrt{c_2(w_1^2-3 w_1+1.5)}\right)^2+\frac{w_1}{3(N-1)}\leq\frac{c_1N_1 + c_2N_2}{4C(N-1)} . \eqno(2.30)$$
We have
$$\eqalign{
\frac{N_2}{6 C  (N_2-1)}&\left(w_1\sqrt{c_1}+\sqrt{c_2(w_1^2-3 w_1+1.5)}\right)^2+\frac{w_1}{3(N-1)}\leq\cr
&\frac{1}{6C(N-1)}\left(N\left(w_1\sqrt{c_1}+\sqrt{c_2(w_1^2-3 w_1+1.5)}\right)^2+2w_1C\right) . \cr
} \eqno(2.31)$$
Now
$$\eqalign{
&\frac{c_1N_1 + c_2N_2}{4C(N-1)}-\frac{1}{6C(N-1)}\left(N\left(w_1\sqrt{c_1}+\sqrt{c_2(w_1^2-3 w_1+1.5)}\right)^2+2w_1C\right)=\cr
&\frac{w_1}{12C(N-1)}\left(N\left((c_1+c_2)(3-2w_1)-4\sqrt{c_1c_2((w_1-3)w_1+1.5)}\right)-4C\right) . \cr
}\eqno(2.32)$$
Since for $0\leq w_1\leq w_1^*$

\item{1.} $(c_1+c_2)(3-2w_1)>0$ is decreasing in $w_1$

\item{2.} $4\sqrt{c_1c_2((w_1-3)w_1+1.5)}>0$ is decreasing in $w_1$

\item{3.} $(c_1+c_2)(3-2w_1)>4\sqrt{c_1c_2((w_1-3)w_1+1.5)}$ for all $w_1$

\item{4.} for $w_1=w_1^*$:  $(c_1+c_2)(3-2w_1)-4\sqrt{c_1c_2((w_1-3)w_1+1.5)}\geq4C$

we have for all $0\leq w_1\leq w_1^*$
$$\frac{w_1}{12C(N-1)}\left(N\left((c_1+c_2)(3-2w_1)-4\sqrt{c_1c_2((w_1-3)w_1+1.5)}\right)-4C\right)\geq0 \eqno(2.33)$$
and hence $D_{0.5}^2\tw\leq D_{0.5}^2\tc$.

For $w_1> w_1^*$ the maximal variance of the estimator $\tw$ equals
%
%
%
%
%
%
%
%
%
%
%
%
%
%
%
%
%
%
%
$$D_{\theta^*}^2\tw=\frac{\big(n_1 (c_1 b_2+N_2c_2  (3N_2-1) (N_1-1))-C b_2\big)^2}{48 n_1 N^2 (N_2-1) (N_1-1) (C-c_1 n_1) (n_1 (c_1 b_1+c_2 N_2 (N_1-1))-C b_1))}, \eqno(2.34)$$

where

$b_1 = n_1 (2-2 N (2-3 N_2
   w_1))+N_1 (1-N_2) (3 N_1-1)$ and $b_2= n_1 (N-2)-N_1 (N_2-1)$.

Similar calculations as in the case $w_1\leq w_1^*$ show that $D_{\theta^*}^2\tw\leq D_{0.5}^2\tc$.
\bigskip

\section{3. Numerical results}
Table 2 shows certain numerical results for $N= 30000$, $c_1=1$, $c_2=3$ and $C=1200$.
$$\vbox{\tabskip1em \offinterlineskip \halign{
\strut\hfil$#$\hfil&&#\vrule&\hfil$#$\hfil\cr
\multispan{15}{\bf Table 2.} Maximal variances $D_{\widetilde{\theta}}^2 \tw, C=1200,c_1=1,c_2=3, N=30000$\hfil\cr\noalign{\vskip5pt}
w_1 && n_1^{opt} &&n_w &&n_c && D_{\widetilde{\theta}}^2 \tw  && D^2_{0.5} \tc &&reduction \cr \noalign{\hrule}
0.05&&\phantom{0}29&&419&&413&&0.0005837&&0.0005970&&\phantom{0}2.23\%\cr
0.10&&\phantom{0}59&&439&&428&&0.0005504&&0.0005757&&\phantom{0}4.40\%\cr
0.15&&\phantom{0}92&&461&&444&&0.0005169&&0.0005547&&\phantom{0}6.82\%\cr
0.20&&127&&484&&461&&0.0004828&&0.0005339&&\phantom{0}9.57\%\cr
0.25&&165&&510&&480&&0.0004488&&0.0005129&&12.54\%\cr
0.30&&207&&538&&500&&0.0004127&&0.0004916&&16.05\%\cr
0.35&&252&&568&&521&&0.0003767&&0.0004712&&20.23\%\cr
0.40&&327&&618&&545&&0.0003056&&0.0004503&&32.15\%\cr
0.45&&383&&655&&571&&0.0002984&&0.0004295&&30.53\%\cr
0.50&&439&&692&&600&&0.0002853&&0.0004083&&30.13\%\cr
}}$$
In the first column of Table 2 values of $w_1$ are given. In the second column the optimal number of units from the first strata in the sample is shown. It is a value $n_1^{opt}$, which gives minimum of $D_{\widetilde{\theta}}^2\tw$ (an exemplary numerical code in Mathematica for calculating optimal allocation is given in Appendix) . Column $n_w$ shows the total sample size: $n_1^{opt}+n_2^{opt}$.  The values of $n_c$~are given in the fourth column. The next column contains minimal (maximal) variance $D_{\widetilde{\theta}}^2 \tw$. In the one before the last column the values of maximal variances $D^2_{0.5} \tc$ are given. In the last column
$$
reduction = \left( 1 - {D_{\widetilde{\theta}}^2 \tw\over D^2_{0.5} \tc} \right)\cdot 100\% \eqno(3.1)
$$
is given.
In considered numerical example, for each value of $w_1$, the maximal variance of the estimator $\tc$ is greater than the maximal variance of the estimator $\tw$ with averaged sample allocation. Furthermore, total sample sizes for stratified random sampling are not smaller than for simple random sampling.
In Figures 1 and 2 variance of $\tw$, as well as the variance of $\tc$ are drawn for $w_1>w_1^*$ and $w_1\leq w_1^*$, respectively and for optimal $n_1$.

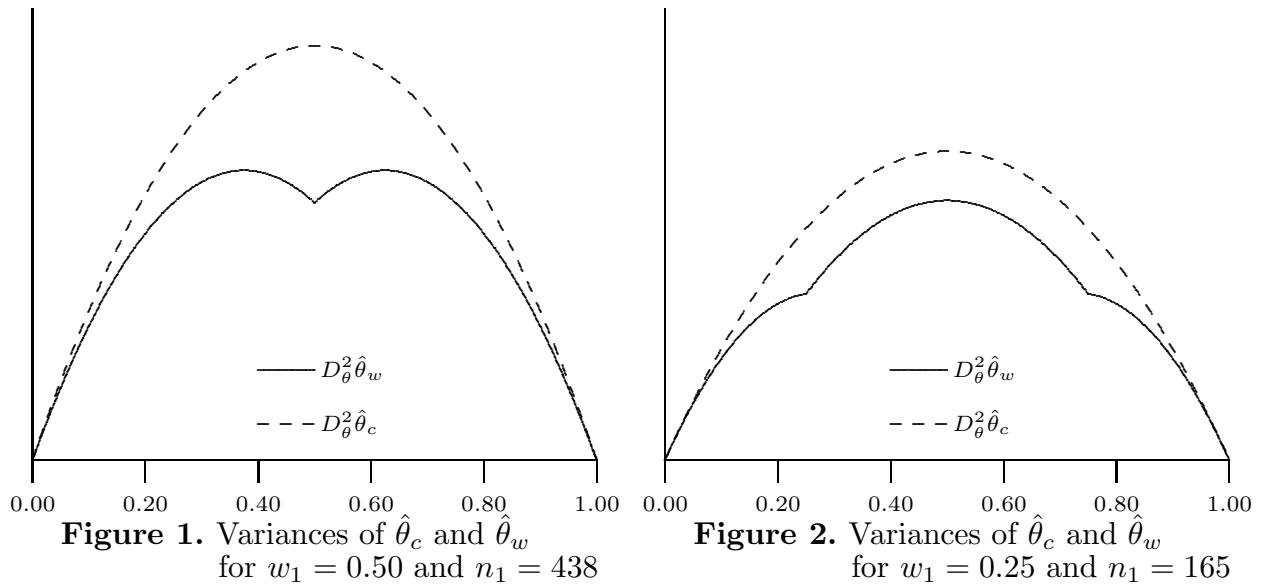
\setbox102=\hbox{%
\beginpicture
\sevenrm
\setcoordinatesystem units <75truemm,0.24truecm>
\setplotarea x from 0 to 1, y from 0.00 to 25
\axis bottom
ticks numbered from 0 to 1 by 0.20 /
\axis left /

\plot "n1-438.tex" \plot 0.4 5 0.5 5 /
 \put{$\scriptstyle D_{\theta}^2\tw$} [l] at 0.51 5
 \setdashpattern<1.5truemm,1.5truemm>
 \plot "calaw50.tex"
 \plot 0.4 2 0.5 2 /
 \put{$\scriptstyle D_{\theta}^2\tc$} [l] at 0.51 2
 \put {{\bf Figure 1.} \rm Variances of $\tc$ and $\tw$} [l] at 0.05 -4
 \put {\rm for $w_1 = 0.50$ and $n_1=438$ } [l] at 0.33 -6
\endpicture
}

\setbox103=\hbox{%
\beginpicture
\sevenrm
\setcoordinatesystem units <75truemm,0.24truecm>
\setplotarea x from 0 to 1, y from 0.00 to 25
\axis bottom
ticks numbered from 0 to 1 by 0.20 /
\axis left /

\plot "n1-180.tex"
\plot 0.4 5 0.5 5 /
\put{$\scriptstyle D_{\theta}^2\tw$} [l] at 0.51 5
\setdashpattern<1.5truemm,1.5truemm>
 \plot "cala.tex"
 \plot 0.4 2 0.5 2 /
 \put{$\scriptstyle D_{\theta}^2\tc$} [l] at 0.51 2
 \put {{\bf Figure 2.} \rm Variances of $\tc$ and $\tw$} [l] at 0.05 -4
 \put {\rm for $w_1 = 0.25$ and $n_1=165$ } [l] at 0.33 -6
\endpicture
}


\line{\copy102\hss\copy103}
\bigskip

\section{4. Summary}
In the paper some approach to optimal sample allocation with respect to limited funds was proposed. Two estimators of an unknown fraction $\theta$ in the finite population were considered: the standard estimator $\tc$ and the stratified estimator $\tw$. It was shown that both estimators are unbiased. Their variances were compared. It was proved that 'the worst' variance of $\tw$ with proposed sample allocation is smaller than 'the worst' variance of $\tc$. The numerical example was presented. In that example it was shown that 'the worst' variance of the stratified estimator may be smaller up to $30\%$ than 'the worst' variance of the classical estimator. For such approach there is no need to estimate unknown $\theta_1$ and $\theta_2$ by preliminary sample. 

\section{\bf References}

\everypar = {
\parindent=0pt
\hangindent=8mm
\hangafter=1
}\noindent

\def\textit#1{{\it#1}}

\noindent  Armitage, P., {\it A Comparison of Stratified with Unrestricted Random Sampling from a~Finite Population}, Biometrika, 34, 3/4 , 273-280, {\bf1947}.


\noindent  Cochran, W. G.,  {\it Sampling Techniques (3rd ed.)}, New York: John Wiley, {\bf 1977}.
\vskip 3mm

\noindent  Hansen, M. H., Hurwitz, W. N., {\it The problem of non-response in sample surveys}, Journal of the American Statistical Association, 41, 236, 517-529, {\bf1946}.

\noindent  Horgan, J. M., {\it Stratification of Skewed Populations: A review}, International Statistical Review, 74, 1, 67-76, {\bf2006}.

\noindent  Johnson, N. L. Kotz, S., {\it Discrete distributions: distributions in statistics}, Houghton Mifflin Company, Boston, {\bf 1969}.
\vskip 3mm
 \noindent  Neyman, J., {\it On the two different aspects of the representative method: The method of stratified sampling and the method of purposive selection}, Journal of the Royal Statistical Society, 97, 558-606, {\bf1934}.


\noindent  Sieradzki, D. Zieli{\'n}ski, W., {\it Sample allocation in estimation of proportion in finite population divided among two strata}, Statistics in Transition new series, 18, 3, 541-548, DOI 10.21307, {\bf2017}.


\noindent  Zieli{\'n}ski, W., {\it A remark on estimating defectiveness in sampling acceptance inspection}, Colloquium Biometricum, 46, 9-14, {\bf2016}.

\section{Appendix}
\everypar = {
\parindent=0pt
\hangindent=0mm
\hangafter=4
}\noindent

\vskip 3mm
An exemplary Mathematica code for calculating optimal allocation is enclosed. Naturally, one can also use other mathematical or
statistical packages (in a similar way) to find the value of $n_1^{opt}$.
\vskip 3mm
\baselineskip=1.5pt

\font\ttmale=pltt8

\begingroup
\parindent0pt
\obeylines
\ttmale

In[1]:=
(*input*)
M = 30000;(*population size*)
K = 1200;(*available funds*)
w1 = 0.25;(*first strata weight*)
k1 = 1; k2 = 3;(*costs of sampling*)
(*end of input*)

M1 = w1*M; M2 = M - M1;

(*variance for 0<$\scriptstyle\backslash$[Theta]<w1*)
VarianceLeft[n1$\_$,n2$\_$,$\scriptstyle\backslash$[Theta]$\_$]=$\scriptstyle\backslash$[Theta]/(6M)
\hskip0.5em(((3M1-1)(M1-n1))/((M1-1)n1)+
\hskip0.5em((3M2-1)(M2-n2))/((M2-1)n2)-2$\scriptstyle\backslash$[Theta]M((M1-n1)/((M1-1)n1)+(M2-n2)/((M2-1)n2)));

(*variance for w1<$\scriptstyle\backslash$[Theta]<1-w1*)
VarianceMiddle[n1$\_$,n2$\_$,$\scriptstyle\backslash$[Theta]$\_$]=
\hskip0.5em(2n1n2w1(M+1)+M1M2((n1+n2)w1-3n1)-M1(n2w1+n1(1-w1)))/
\hskip0.5em(6Mn1n2(M2-1))+((M2-n2)(1-$\scriptstyle\backslash$[Theta])$\scriptstyle\backslash$[Theta])/(n2(M2-1));

(*variance for 0<$\scriptstyle\backslash$[Theta]<1*)
VarianceStratified[n1$\_$,$\scriptstyle\backslash$[Theta]$\_$]=
\hskip0.5emIf[0<$\scriptstyle\backslash$[Theta]<w1,VarianceLeft[n1,(K-k1*n1)/k2,$\scriptstyle\backslash$[Theta]],
\hskip1emIf[w1<=$\scriptstyle\backslash$[Theta]<=1-w1,VarianceMiddle[n1,(K-k1*n1)/k2,$\scriptstyle\backslash$[Theta]],
\hskip1.5emVarianceLeft[n1,(K-k1*n1)/k2,1-$\scriptstyle\backslash$[Theta]]]];

(*determining optimal n1*)
a=1;b=K/k1;
coefficient=N[(Sqrt[5]-1)/2];
xL=b-coefficient*(b-a);
xR=a+coefficient*(b-a);
eps=0.1;
While[(b-a)>eps,
$\scriptstyle\{$LL=FindMaximum[VarianceStratified[xL,$\scriptstyle\backslash$[Theta]],$\scriptstyle\{$$\scriptstyle\backslash$[Theta],0.1$\scriptstyle\}$][[1]];
RR=FindMaximum[VarianceStratified[xR,$\scriptstyle\backslash$[Theta]],$\scriptstyle\{$$\scriptstyle\backslash$[Theta],0.1$\scriptstyle\}$][[1]];
If[LL<RR,
$\scriptstyle\{$b=xR;xR=xL;xL=b-coefficient*(b-a);$\scriptstyle\}$,
$\scriptstyle\{$a=xL;xL=xR;xR=a+coefficient*(b-a);$\scriptstyle\}$];
$\scriptstyle\}$];
n1opt=IntegerPart[(a+b)/2];

(*output*)
v1=FindMaximum[VarianceStratified[n1opt,$\scriptstyle\backslash$[Theta]],$\scriptstyle\{$$\scriptstyle\backslash$[Theta],0.1$\scriptstyle\}$][[1]];
nsdwr=IntegerPart[K/(w1*k1+(1-w1)*k2)];
v2=N[(M-nsdwr)/(4*(M-1)*nsdwr)];
Print["optimal n1: ", n1opt]
Print["optimal n2: ", IntegerPart[(K - k1*n1opt)/k2]]
Print["maximal stratified variance: ", v1]
Print["simple draw without replacement sample size: ", nsdwr]
Print["simple draw without replacement maximal variance: ", v2]
Print["variance reduction: ", (1 - v1/v2)*100 "\%"]
\medskip
Out[1]:=
optimal n1: 165
optimal n2: 345
maximal stratified variance: 0.000448249
simple draw without replacement sample size: 480
simple draw without replacement maximal variance: 0.000512517
variance reduction: 12.5397\%

\endgroup

\bye